\documentclass{article}
\usepackage{amsmath,amssymb,amsthm}

\newcommand{\N}{{\mathbb N}}

\newcommand{\Z}{{\mathbb Z}}
\newcommand{\C}{{\mathbb C}}

\makeatletter
\def\eqnarray{\stepcounter{equation}\let\@currentlabel=\theequation
\global\@eqnswtrue
\tabskip\@centering\let\\=\@eqncr
$$\halign to \displaywidth\bgroup\hfil\global\@eqcnt\z@
  $\displaystyle\tabskip\z@{##}$&\global\@eqcnt\@ne
  \hfil$\displaystyle{{}##{}}$\hfil
  &\global\@eqcnt\tw@ $\displaystyle{##}$\hfil
  \tabskip\@centering&\llap{##}\tabskip\z@\cr}

\def\endeqnarray{\@@eqncr\egroup
      \global\advance\c@equation\m@ne$$\global\@ignoretrue}

\def\@yeqncr{\@ifnextchar [{\@xeqncr}{\@xeqncr[5pt]}}
\makeatother

\begin{document}
\bibliographystyle{plain}

\centerline{\LARGE{\bf Derivatives with respect to the degree}}
\vskip 0.1 truecm
\centerline{\LARGE{\bf and order of associated Legendre functions}}
\vskip 0.1 truecm
\centerline{\LARGE{\bf for $|z|>1$ using modified Bessel functions}}
\vskip 0.5 truecm
\centerline{\Large{Howard S. Cohl}}
\vskip 0.4 truecm
\centerline{\normalsize{Department of Mathematics, University of Auckland, 38 Princes Street, 
Auckland, New Zealand}}

\vspace{0.3cm}

\begin{center}
(Received: 09 September 2009; Revised [Final]: 25 October 2009)
\end{center}
\vspace{0.0cm}

\begin{center}
\begin{minipage}{0.7\textwidth}
Expressions for the derivatives with respect to order of modified
Bessel functions evaluated at integer orders and certain integral 
representations of associated Legendre functions with modulus argument 
greater than unity are used to compute derivatives of the associated Legendre
functions with respect to their parameters.  For the associated Legendre 
functions of the first and second kind, derivatives with respect to the 
degree are evaluated at odd-half-integer degrees, for general complex 
orders, and derivatives with respect to the order are evaluated at 
integer orders, for general complex degrees.
\end{minipage}
\end{center}
\vspace{-0.2cm}
\begin{center}
\begin{minipage}{0.7\textwidth}
\end{minipage}
\end{center}
\vspace{-0.2cm}
\begin{center}
\begin{minipage}{0.7\textwidth}
\end{minipage}
\end{center}

\vspace{0.2cm}

$\ ^\ast$Email:~h.cohl@math.auckland.ac.nz

\newpage

\section{Introduction}

\noindent Here we present formulae for derivatives of associated Legendre 
functions (hereafter referred to as Legendre functions) of the first kind $P_\nu^\mu(z)$ and 
the second kind $Q_\nu^\mu(z)$, with respect to their parameters, namely the degree $\nu$ and the 
order $\mu$.
Some formulae relating to these derivatives have been previously 
noted \cite{MOS}
and also there has been recent work in this area 
\cite{Brych,Szmy1,Szmy2,Szmy3,Szmy4}
with Brychkov (2009) \cite{Brychkov} giving a recent reference covering the
regime for argument 
%$|z|<1$. 
$z\in[-1,1]$. 
In this paper we cover parameter derivatives of Legendre
functions for argument $|z|>1.$ 

The strategy applied in this paper is to incorporate 
derivatives with respect to order evaluated at integer orders
for modified Bessel functions
(see \cite{BryGed,Abra,MOS})
to compute derivatives with respect to the degree and the order of Legendre functions.
Below we apply these results through certain integral representations of 
Legendre functions in terms of modified Bessel functions.

\section{A useful transformation on the complex plane}

There is a transformation over an open subset of the complex plane which is particularly
useful in studying Legendre functions \cite{Hob,Abra}.
This transformation, which is valid on a certain domain of the complex numbers,
accomplishes the following
\begin{eqnarray}
&\cosh z\leftrightarrow\coth w&\nonumber\\
&\coth z\leftrightarrow\cosh w&\nonumber\\
&\sinh z\leftrightarrow 1/\sinh w,&\nonumber
\end{eqnarray}
\noindent where $\cosh$, $\coth,$ and $\sinh$ are the complex hyperbolic
cosine, cotangent, and sine functions respectively.
This transformation is accomplished using the following map
\begin{equation}
w(z)=\log\coth\frac{z}{2},
\label{map}
\end{equation}
(where $\log$ is the complex natural logarithm) which is verified to be an involution over an open subset of the complex plane given by
$-\pi < \mbox{Im}(z) < \pi$, where one removes a branch given by $\mbox{Re}(z)\le 0$, 
$\mbox{Im}(z)=0.$ This mapping is $\pi$ periodic in the imaginary direction and is 
locally injective over the entire complex plane when the following set of complex
numbers are removed $\{z: z=i\pi n$, $n\in\Z\} \cup \{z: \mbox{Im}(z)=2\pi n \mbox{\ and\ }
\mbox{Re}(z)\le 0, \ n\in\Z, z\in\C\}$.

This transformation is particularly useful for certain Legendre functions
which have natural domain given by the real interval $(1,\infty)$, such as toroidal 
harmonics \cite{CT,CRTB} (and for other Legendre functions which one might 
encounter in potential theory), Legendre functions of 
the first and second kind with odd-half-integer degree and integer order.  The real argument 
of these Legendre functions naturally occur in $[1,\infty)$, and these are the simultaneous ranges 
of both the real hyperbolic cosine and cotangent functions.  
%Other applications 
%would be general Legendre functions with arguments given by the natural domain 
%of this transformation.  

One application of this map occurs with the Whipple 
transformation of Legendre functions \cite{Whip,CTRS}
under index (degree and order) interchange.  See for instance, eqs. (8.2.7) and (8.2.8) 
in \cite{Abra}, namely

\begin{equation}
P_{-\mu-1/2}^{-\nu-1/2}
\biggl(\frac{z}{\sqrt{z^2-1}}\biggr)=
\sqrt{\frac{2}{\pi}}
\frac{(z^2-1)^{1/4}\mathrm{e}^{-i\mu\pi}}
{\Gamma(\nu+\mu+1)}Q_\nu^\mu(z),
\label{whipple}
\end{equation}

%%\noindent and,
%%\begin{equation}
%%Q_{-\mu-1/2}^{-\nu-1/2}
%%\biggl(\frac{z}{\sqrt{z^2-1}}\biggr)=
%%-i (\pi/2)^{1/2} \Gamma(-\nu-\mu) (z^2-1)^{1/4}\mathrm{e}^{-i\nu\pi} 
%%P_\nu^\mu(z),
%%\]
\noindent 
which are valid for $\mathrm{Re}\ \!z>0$ 
and for all complex $\nu$ and $\mu,$ except where the
functions are not defined. $\Gamma(z)$ is the Gamma function (see \cite{Abra}).

\section{Parameter derivative formulas from $K_\nu(t)$}

Starting with Gradshteyn \& Ryzhik (2007) eq.~(6.628.7) \cite{Grad} we have
\begin{eqnarray}
\displaystyle \int_0^\infty e^{-zt}K_\nu(t)t^{\alpha-1/2}dt&=&\sqrt{\frac{\pi}{2}}
\Gamma\left(\alpha-\nu+\frac12\right)\Gamma\left(\alpha+\nu+\frac12\right)
\left(z^2-1\right)^{-\alpha/2}
P_{\nu-1/2}^{-\alpha}(z)\nonumber\\[0.2cm]
&=&\Gamma\left(\alpha-\nu+\frac12\right)\left(z^2-1\right)^{-\alpha/2-1/4}e^{-i\pi\nu}
Q_{\alpha-1/2}^\nu\left(\frac{z}{\sqrt{z^2-1}}\right),
\label{IK}
\end{eqnarray}
where $K_\nu(t)$ is a modified Bessel function of the second kind with 
order $\nu$, and 
the two equalities are established through the Whipple transformation,
eq.~(\ref{whipple}).  

We would like to generate an analytical expression for the
derivative of the Legendre function of the second kind with respect to its
order, evaluated at integer orders.  In order to do this our strategy is to 
solve the above integral expression for the Legendre function of the second kind, 
differentiate with respect to the order, evaluate at integer orders, and take 
advantage of the corresponding formula for differentiation with respect to order 
for modified Bessel functions of the second kind (see \cite{Brychkov,BryGed,Abra,MOS}).
Using the expression for the Legendre function of the second kind in 
eq.~(\ref{IK}), we solve for $Q_{\nu-1/2}^\mu(z)$ and re-express using the map in
eq.~(\ref{map}).  This gives us the following expression

\[
Q_{\nu-1/2}^\mu(z)=\frac{\left(z^2-1\right)^{-\nu/2-1/4}e^{i\pi\mu}}{\Gamma\left(\nu-\mu+\frac12\right)}
\int_0^\infty \exp\left(\frac{-zt}{z^2-1}\right)K_\mu(t)t^{\nu-1/2}dt.
\]

\noindent Differentiating with respect to the order $\mu$ and evaluating at 
$\mu=\pm m$, where $m\in\N_0=\N\cup\{0\}$ yields
\begin{eqnarray}
&\displaystyle \left[\frac{\partial}{\partial \mu}Q_{\nu-1/2}^\mu(z)\right]_{\mu=\pm m}=
\left(z^2-1\right)^{-\nu/2-1/4}
\left[
\frac{\partial}{\partial\mu}
\frac{e^{i\pi\mu}}{\Gamma\left(\nu-\mu+\frac12\right)}
\right]_{\mu=\pm m}
\int_0^\infty \exp\left(\frac{-zt}{\sqrt{z^2-1}}\right)K_{\pm m}(t)t^{\nu-1/2}dt&
\nonumber\\[0.2cm]
&\displaystyle +
\frac{\left(z^2-1\right)^{-\nu/2-1/4}(-1)^m}{\Gamma\left(\nu\mp m+\frac12\right)}
\int_0^\infty \exp\left(\frac{-zt}{\sqrt{z^2-1}}\right)t^{\nu-1/2}
\left[
\frac{\partial}{\partial\mu}
K_\mu(t)
\right]_{\mu=\pm m}dt.&
\nonumber
\end{eqnarray}

\noindent The derivative from the first term is given as
\[
\left[
\frac{\partial}{\partial\mu}
\frac{e^{i\pi\mu}}{\Gamma\left(\nu-\mu+\frac12\right)}
\right]_{\mu=\pm m}=
\frac{(-1)^m}{\Gamma\left(\nu\mp m+\frac12\right)}
\left[i\pi+\psi\left(\nu\mp m+\frac12\right) \right],
\]
\noindent where the Digamma function $\psi(z)$ is defined in terms of the derivative of
the Gamma function with respect to its argument $z$ through
\[
\Gamma^\prime(z)=\Gamma(z)\psi(z).
\]
\noindent The derivative in the second integral 
(see \cite{Brychkov,BryGed,Abra,MOS})
is given by
\begin{equation}
\left[
\frac{\partial}{\partial\mu}
K_\mu(t)
\right]_{\mu=\pm m}=\pm m!
\sum_{k=0}^{m-1}
\frac{1}{k!(m-k)}\frac{t^{k-m}}{2^{k-m+1}}K_k(t) 
\label{dKdn}
\end{equation}
\noindent 
(see for instance eq.~(1.14.2.2) in Brychkov (2008) \cite{Brychkov}).
Substituting these expressions for the derivatives into the two
integrals and using the map in eq.~(\ref{map}) to re-evaluate 
these integrals in terms of Legendre functions gives 
the following general expression for the derivative of the 
Legendre function of the second kind with respect to its order evaluated at
integer orders as
\begin{eqnarray}
&\displaystyle 
\frac{\Gamma(\nu\mp m+\frac12)}{\Gamma(\nu-m+\frac12)}\left[\frac{\partial}{\partial \mu}Q_{\nu-1/2}^\mu(z)\right]_{\mu=\pm m}=
\left[i\pi+\psi\left(\nu\mp m+\frac12\right) \right] Q_{\nu-1/2}^m(z)&\nonumber \\[0.2cm] 
&\displaystyle \pm m!\sum_{k=0}^{m-1}\frac{(-1)^{k-m}\left(z^2-1\right)^{(k-m)/2}}{k!(m-k)2^{k-m+1}}
Q_{\nu+k-m-1/2}^k(z).\nonumber
\end{eqnarray}
\noindent For $\mu=0$ there is no contribution from the sum and the result is 
\[
\left[\frac{\partial}{\partial \mu}Q_{\nu-1/2}^\mu(z)\right]_{\mu=0}=
\left[ i\pi+\psi\left(\nu+\frac12\right) \right] Q_{\nu-1/2}(z),
\]
\noindent which agrees with that given in Magnus, Oberhettinger,
and Soni (1966) \cite{MOS} in \S4.4.3.  
We are now able to obtain formulas for non-zero values of $\mu$ such as for $\mu=-1$
\[
\left(\nu^2-\frac14\right)
\left[\frac{\partial}{\partial \mu}Q_{\nu-1/2}^\mu(z)\right]_{\mu=-1}=
\left[i\pi+\psi\left(\nu+\frac32\right) \right] Q_{\nu-1/2}^1(z)
+\left(z^2-1\right)^{-1/2}Q_{\nu-3/2}(z),
\]
\noindent or for $\mu=+1$
\[
\left[\frac{\partial}{\partial \mu}Q_{\nu-1/2}^\mu(z)\right]_{\mu=1}=
\left[i\pi+\psi\left(\nu-\frac12\right) \right] Q_{\nu-1/2}^1(z)
-\left(z^2-1\right)^{-1/2}Q_{\nu-3/2}(z),
\]
\noindent or for other integer values of $\mu$.

If we start with the expression for the Legendre function of
the first kind in eq.~(\ref{IK}) and solve for $P_{\nu-1/2}^{-\mu}(z)$ we have

\begin{equation}
P_{\nu-1/2}^{-\mu}(z)=\sqrt{\frac{2}{\pi}}
\frac{\left(z^2-1\right)^{\mu/2}}
{\Gamma\left(\mu-\nu+\frac12\right)\Gamma\left(\mu+\nu+\frac12\right)}
\int_0^\infty e^{-zt}K_\nu(t)t^{\mu-1/2}dt.
\label{PK}
\end{equation}

\noindent Differentiating with respect to degree $\nu$ and evaluating at 
$\nu=\pm n$, where $n\in\N_0$ yields
\begin{eqnarray}
\left[\frac{\partial}{\partial \nu}P_{\nu-1/2}^{-\mu}
(p)\right]_{\nu=\pm n}=
\sqrt{\frac{2}{\pi}}\left(z^2-1\right)^{\mu/2}
\left[\frac{\partial}{\partial \nu}\frac{1}
{\Gamma\left(\mu-\nu+\frac12\right)\Gamma\left(\mu+\nu+\frac12\right)}
\right]_{\nu=\pm n}
\int_0^\infty e^{-zt}K_{\pm n}(t)t^{\mu-1/2}dt\nonumber\\[0.2cm]
+
\sqrt{\frac{2}{\pi}}
\frac{\left(z^2-1\right)^{\mu/2}}
{\Gamma\left(\mu\mp n+\frac12\right)\Gamma\left(\mu\pm n+\frac12\right)}
\int_0^\infty e^{-zt}t^{\mu-1/2}
\left[\frac{\partial}{\partial \nu}K_\nu(t)\right]_{\nu=\pm n}dt.
\nonumber
\end{eqnarray}

\noindent The derivative from the first term is given as
\[
\left[
\frac{\partial}{\partial\nu}
\frac{1}{\Gamma\left(\mu-\nu+\frac12\right)\Gamma\left(\mu+\nu+\frac12\right)}
\right]_{\nu=\pm n}=
\frac{\psi\left(\mu\mp n+\frac12\right)-\psi\left(\mu\pm n+\frac12\right)}
{\Gamma\left(\mu\pm n+\frac12\right)\Gamma\left(\mu\mp n+\frac12\right)}.
\]

\noindent Substituting this expression for the derivative and that given
in eq.~(\ref{dKdn}) yields the following general expression for the 
derivative of the Legendre function of the first kind with 
respect to its degree evaluated at odd-half-integer degrees as

\begin{eqnarray}
\pm\left[\frac{\partial}{\partial \nu}P_{\nu-1/2}^{-\mu}
(z)\right]_{\nu=\pm n}=
\left[
\psi\left(\mu-n+\frac12\right)-\psi\left(\mu+n+\frac12\right)
\right]
P_{n-1/2}^{-\mu}(z)\nonumber\\[0.2cm]
+\frac{n!}{\Gamma\left(\mu+n+\frac12\right)}
\sum_{k=0}^{n-1}
\frac{\Gamma\left(\mu-n+2k+\frac12\right)\left(z^2-1\right)^{(n-k)/2}}
{k!(n-k)2^{k-n+1}}
P_{k-1/2}^{-\mu+n-k}(z).\nonumber
\end{eqnarray}

\noindent If one makes a global replacement, $-\mu\mapsto\mu$, using
the properties of Gamma and Digamma functions, this result reduces to
\begin{eqnarray}
\pm\left[\frac{\partial}{\partial \nu}P_{\nu-1/2}^{\mu}
(z)\right]_{\nu=\pm n}=
\left[
\psi\left(\mu+n+\frac12\right)-\psi\left(\mu-n+\frac12\right)
\right]
P_{n-1/2}^{\mu}(z)\nonumber\\[0.2cm]
+n!\ \Gamma\left(\mu-n+\frac12\right)
\sum_{k=0}^{n-1}
\frac{\left(z^2-1\right)^{(n-k)/2}}
{\Gamma\left(\mu+n-2k+\frac12\right)k!(n-k)2^{k-n+1}}
P_{k-1/2}^{\mu+n-k}(z).\nonumber
\end{eqnarray}

\noindent Note that by using the recurrence relation for Digamma functions (see \S 1.2 in \cite{MOS}) 
\[
\psi(z+1)=\psi(z)+\frac{1}{z},
\]
we can establish
\[
\psi\left(\mu+n+\frac12\right)-\psi\left(\mu-n+\frac12\right)=
2\mu\sum_{l=1}^n\left[\mu^2-\left(l-\frac12\right)^2\right]^{-1}.
\]

\noindent For $\nu=0$ there is no contribution from the sum and the result is 
\[
\left[\frac{\partial}{\partial \nu}P_{\nu-1/2}^{\mu}
(z)\right]_{\nu=0}=0,
\]
\noindent which agrees with that given in Magnus, Oberhettinger,
and Soni (1966) \cite{MOS} in \S4.4.3.  
We are now able to obtain formulas for non-zero values of $\nu$ such as for $\nu=\pm 1$
\[
\pm\left(\mu^2-\frac14\right)
\left[\frac{\partial}{\partial \nu}P_{\nu-1/2}^{\mu}
(z)\right]_{\nu=\pm 1}=
2\mu P_{1/2}^{\mu}(z)+
\left(z^2-1\right)^{1/2}
P_{-1/2}^{\mu+1}(z),
\]
\noindent or for other integer values of $\nu$.

Note that this method does not seem amenable to
computing derivatives with respect to the degree of Legendre functions of
the form $P_{\nu}^{\mu}(z)$ evaluated at integer degrees, since shifting the degree by $+1/2$ in 
eq.~(\ref{PK})
converts the modified Bessel function of the second kind to a form like
$K_{\nu+1/2}(t)$ and the derivative with respect to order of this
Bessel function (see \cite{BryGed,Abra,MOS})
is not of a form which is easily integrated.

\section{Parameter derivative formulas from $I_\nu(t)$}

Starting with Gradshteyn \& Ryzhik (2007) eq.~(6.624.5) \cite{Grad}, see also
Prudnikov et.~al.~(1988) \cite{Prud}, we have
\begin{eqnarray}
\displaystyle \int_0^\infty e^{-zt}I_\nu(t)t^{\alpha-1/2}dt&=&
\sqrt{\frac{2}{\pi}}
e^{-i\pi\alpha}\left(z^2-1\right)^{-\alpha/2}
Q_{\nu-1/2}^\alpha(z)\nonumber\\[0.2cm]
&=&\Gamma\left(\alpha+\nu+\frac12\right)\left(z^2-1\right)^{-\alpha/2-1/4}
P_{\alpha-1/2}^{-\nu}\left(\frac{z}{\sqrt{z^2-1}}\right),
\label{II}
\end{eqnarray}
where $I_\nu(t)$ is a modified Bessel function of the first kind with order $\nu$, and 
the two equalities are established through the Whipple transformation,
eq.~(\ref{whipple}).  We will use this particular integral representation of
Legendre functions to compute certain derivatives of 
the Legendre functions with respect to the degree and order.

We start with the integral representation of the Legendre
function of the second kind in eq.~(\ref{II}).
Differentiating with respect to the degree $\nu$ and evaluating at $\nu=\pm n$, where 
$n\in\N_0$, one obtains
\begin{equation}
\left[\frac{\partial}{\partial \nu} Q_{\nu-1/2}^\mu(z)
\right]_{\nu=\pm n}
=\sqrt{\frac{\pi}{2}}e^{i\pi\mu}\left(z^2-1\right)^{\mu/2}
\int_0^\infty e^{-zt} t^{\mu-1/2}
\left[
\frac{\partial}{\partial \nu} I_\nu(t)\right]_{\nu=\pm n}dt.
\label{dQn}
\end{equation}

\noindent The derivative of the modified Bessel function of the first kind
in eq.~(\ref{dQn}) (see \cite{Brychkov,BryGed,Abra,MOS}) is given by
\begin{equation}
\left[\frac{\partial}{\partial \nu} I_\nu(t)\right]_{\nu=\pm n}=
(-1)^{n+1}K_n(t)\pm n!\sum_{k=0}^{n-1}
\frac{(-1)^{k-n}}{k!(n-k)}
\frac{t^{k-n}}{2^{k-n+1}}
I_k(t)
\label{dIdn}
\end{equation}
\noindent (see for instance eq.~(1.13.2.1) in Brychkov (2008) \cite{Brychkov}).
Inserting eq.~(\ref{dIdn}) into eq.~(\ref{dQn}) and using
eqs.~(\ref{IK}) and (\ref{II}), we obtain the following general expression for the
derivative of the Legendre function of the second kind with respect to its degree
evaluated at odd-half-integer degrees as
\begin{eqnarray}
&\displaystyle
\left[
\frac{\partial}{\partial \nu} Q_{\nu-1/2}^\mu(z)
\right]_{\nu=\pm n}
=-\sqrt{\frac{\pi}{2}}e^{i\pi\mu}
\Gamma\left(\mu-n+\frac12\right)\left(z^2-1\right)^{-1/4}Q_{\mu-1/2}^n
\left(\frac{z}{\sqrt{z^2-1}}\right)&\nonumber\\
&\displaystyle\pm n!
\sum_{k=0}^{n-1}\frac{\left(z^2-1\right)^{(n-k)/2}}{2^{k-n+1}k!(n-k)}
Q_{k-1/2}^{\mu+k-n}(z).&\nonumber
\end{eqnarray}
\noindent For $\nu=0$ there is no contribution from the sum and the result is 
\[
\left[
\frac{\partial}{\partial \nu} Q_{\nu-1/2}^\mu(z)
\right]_{\nu=0}=
-\sqrt{\frac{\pi}{2}}e^{i\pi\mu}\Gamma\left(\mu+\frac12\right)
\left(z^2-1\right)^{-1/4}
Q_{\mu-1/2}\left(\frac{z}{\sqrt{z^2-1}}\right),
\]
which agrees with that given in Magnus, Oberhettinger, and Soni (1966) \cite{MOS} 
in \S 4.4.3.  
We are now able to obtain formulas for non-zero values of $\nu$ such as for $\nu=\pm 1$
\[
\left[
\frac{\partial}{\partial \nu} Q_{\nu-1/2}^\mu(z)
\right]_{\nu=\pm 1}=
-\sqrt{\frac{\pi}{2}}e^{i\pi\mu}\Gamma\left(\mu-\frac12\right)
\left(z^2-1\right)^{-1/4}
Q_{\mu-1/2}^1\left(\frac{z}{\sqrt{z^2-1}}\right)
\pm\left(z^2-1\right)^{1/2}\ Q_{-1/2}^{\mu-1}(z),
\]
\noindent or for other integer values of $\nu$.

We can see that this method does not seem amenable to
computing derivatives with respect to the degree of Legendre functions of
the form $Q_{\nu}^{\mu}(z)$ evaluated at integer degrees, since shifting 
the degree by $+1/2$ in 
eq.~(\ref{dQn})
converts the modified Bessel function of the first kind to a form like
$I_{\nu+1/2}(t)$ and the derivative with respect to order of this
Bessel function (see \cite{BryGed,Abra,MOS})
is not of a form which is easily integrated.

Finally, we obtain a formula for the derivative with
respect to the order for the Legendre function of the first kind evaluated 
at integer orders.
In order to do this we use the integral expression for the 
Legendre function of the first kind given in eq.~(\ref{II}) and
the map given in eq.~(\ref{map}) to convert to the appropriate
argument.  Now use the negative order condition for Legendre 
functions of the first kind (see for example eq.~(22) in \cite{CTRS})
to convert to a positive order.  

Differentiating both sides of the 
resulting expression with respect to the order $\mu$ and evaluating 
at $\mu=\pm m,$ where $m\in\N_0$ yields
\begin{eqnarray}
\displaystyle \left[\frac{\partial}{\partial\mu} P_{\nu-1/2}^\mu(z) \right]_{\mu=\pm m}=
&2&\!\!\hspace{0.3mm}Q_{\nu-1/2}^{\pm m}(z)\nonumber\\[-0.1cm]
\displaystyle &+&\left(z^2-1\right)^{-\nu/2-1/4}
\left\{
\frac{\partial}{\partial\mu}
\left[
\Gamma\left(\nu-\mu+\frac12\right)
\right]^{-1}
\right\}_{\mu=\pm m}\nonumber\\[0.1cm]
&{\,}&\quad\times\int_0^\infty \exp\left(\frac{-zt}{\sqrt{z^2-1}} \right) 
I_{\pm m}(t)t^{\nu-1/2}dt\nonumber\\[0.2cm]
\displaystyle &+&\frac{\left(z^2-1\right)^{-\nu/2-1/4}}{\Gamma\left(\nu\mp m+\frac12\right)}
\int_0^\infty \exp\left( \frac{-zt}{\sqrt{z^2-1}}\right) 
t^{\nu-1/2}
\left[ \frac{\partial}{\partial\mu} I_\mu(t) \right]_{\mu=\pm m}dt.
\nonumber
\end{eqnarray}
\noindent The derivative of the reciprocal of the Gamma function reduces
to $\psi\left(\nu\mp m+1/2\right)/\Gamma\left(\nu\mp m+1/2\right)$.  The derivative
with respect to order for the modified Bessel function of the first kind 
is given in eq.~(\ref{dIdn}).  The integrals are easily obtained by 
applying the map given by eq.~(\ref{map}) as necessary to 
eqs.~(\ref{IK}) and (\ref{II}).  Hence by also using standard properties of 
Legendre, Gamma, and Digamma functions we obtain the
following compact form
\begin{eqnarray}
&\displaystyle \frac{\Gamma(\nu\mp m+\frac12)}{\Gamma(\nu-m+\frac12)}
\left[\frac{\partial}{\partial\mu} P_{\nu-1/2}^\mu(z) \right]_{\mu=\pm m}=
Q_{\nu-1/2}^m(z)+\psi\left(\nu\mp m+\frac12\right)
P_{\nu-1/2}^m(z)&\nonumber\\[0.2cm]
&\displaystyle \pm m!
\sum_{k=0}^{m-1}
\frac{(-1)^{k-m}\left(z^2-1\right)^{(k-m)/2}}{2^{k-m+1}k!(m-k)}
P_{\nu+k-m-1/2}^k(z).&\nonumber
\end{eqnarray}
\noindent For $\mu=0$ there is no contribution from the sum and the result is 
\[
\left[\frac{\partial}{\partial\mu} P_{\nu-1/2}^\mu(z) \right]_{\mu=0}=
Q_{\nu-1/2}(z)+\psi\left(\nu+\frac12\right)P_{\nu-1/2}(z),
\]
\noindent which agrees with that given in Magnus, Oberhettinger,
and Soni (1966) \cite{MOS} in \S4.4.3. 
We are now able to obtain formulas for non-zero values of $\mu$ such as for $\mu=-1$
\[
\left(\nu^2-\frac14\right)
\left[\frac{\partial}{\partial\mu} P_{\nu-1/2}^\mu(z) \right]_{\mu=-1}=
Q_{\nu-1/2}^1(z)+\psi\left(\nu+\frac32\right)P_{\nu-1/2}^1(z)
+\left(z^2-1\right)^{-1/2}P_{\nu-3/2}(z),
\]
\noindent or for $\mu=+1$
\[
\left[\frac{\partial}{\partial\mu} P_{\nu-1/2}^\mu(z) \right]_{\mu=1}=
Q_{\nu-1/2}^1(z)+\psi\left(\nu-\frac12\right)P_{\nu-1/2}^1(z)
-\left(z^2-1\right)^{-1/2}P_{\nu-3/2}(z),
\]
\noindent or for other integer values of $\mu$.

\section{Acknowledgements}
I owe much thanks to the following people: Prof.~Rados{\l}aw Szmytkowski suggested that 
I submit this work for publication and participated in several important discussions. 
I had many valuable discussions with 
Prof.~Ernie Kalnins, who supported me when I was writing this article at the 
University of Waikato. Dr Yury Brychkov engaged me in valuable discussions. 
At the University of Auckland, I had some valuable conversations with 
Dr Garry J.~Tee who carefully proofread this article, and Dr Tom ter Elst 
discussed matters relating to the map in \S2. And I thank the referee for 
helpful remarks.

\bibliography{refbib}

\def\cprime{$'$} \def\cprime{$'$}
\begin{thebibliography}{10}

\bibitem{Abra}
M.~Abramowitz and I.~A. Stegun.
\newblock {\em Handbook of mathematical functions with formulas, graphs, and
  mathematical tables}, volume~55 of {\em National Bureau of Standards Applied
  Mathematics Series}.
\newblock For sale by the Superintendent of Documents, U.S. Government Printing
  Office, Washington, D.C., 1964.

\bibitem{Brych}
Yu.~A. Brychkov.
\newblock {\em {H}andbook of {S}pecial {F}unctions: {D}erivatives, {I}ntegrals,
  {S}eries and {O}ther {F}ormulas}.
\newblock {C}hapman \& {H}all/CRC Press, Boca Raton-London-New York, 2008.

\bibitem{Brychkov}
Yu.~A. Brychkov.
\newblock On the derivatives of the {L}egendre functions ${P}_\nu^\mu(z)$ and
  ${Q}_\nu^\mu(z)$ with respect to $\mu$ and $\nu$.
\newblock {\em Integral Transforms Spec. Funct.},
  21:DOI--10.1080/10652460903069660, 2010.

\bibitem{BryGed}
Yu.~A. Brychkov and K.~O. Geddes.
\newblock On the derivatives of the {B}essel and {S}truve functions with
  respect to the order.
\newblock {\em Integral Transforms Spec. Funct.}, 16(3):187--198, 2005.

\bibitem{CRTB}
H.~S. Cohl, A.~R.~P. Rau, J.~E. Tohline, D.~A. Browne, J.~E. Cazes, and E.~I.
  Barnes.
\newblock Useful alternative to the multipole expansion of $1/r$ potentials.
\newblock {\em Physical Review A: Atomic and Molecular Physics and Dynamics},
  64(5):052509, Oct 2001.

\bibitem{CT}
H.~S. {Cohl} and J.~E. {Tohline}.
\newblock {A Compact Cylindrical Green's Function Expansion for the Solution of
  Potential Problems}.
\newblock {\em The Astrophysical Journal}, 527:86--101, December 1999.

\bibitem{CTRS}
H.~S. {Cohl}, J.~E. {Tohline}, A.~R.~P. {Rau}, and H.~M. {Srivastava}.
\newblock {Developments in determining the gravitational potential using
  toroidal functions}.
\newblock {\em Astronomische Nachrichten}, 321(5/6):363--372, 2000.

\bibitem{Grad}
I.~S. Gradshteyn and I.~M. Ryzhik.
\newblock {\em Table of integrals, series, and products}.
\newblock Elsevier/Academic Press, Amsterdam, seventh edition, 2007.
\newblock Translated from the Russian, Translation edited and with a preface by
  Alan Jeffrey and Daniel Zwillinger, With one CD-ROM (Windows, Macintosh and
  UNIX).

\bibitem{Hob}
E.~W. Hobson.
\newblock {\em The theory of spherical and ellipsoidal harmonics}.
\newblock Chelsea Publishing Company, New York, 1955.

\bibitem{MOS}
W.~Magnus, F.~Oberhettinger, and R.~P. Soni.
\newblock {\em Formulas and theorems for the special functions of mathematical
  physics}.
\newblock Third enlarged edition. Die Grundlehren der mathematischen
  Wissenschaften, Band 52. Springer-Verlag New York, Inc., New York, 1966.

\bibitem{Prud}
A.~P. Prudnikov, Yu.~A. Brychkov, and O.~I. Marichev.
\newblock {\em Integrals and series. {V}ol. 2}.
\newblock Gordon \& Breach Science Publishers, New York, second edition, 1988.
\newblock Special functions, Translated from the Russian by N. M. Queen.

\bibitem{Szmy4}
R.~{Szmytkowski}.
\newblock Addendum to: ``{O}n the derivative of the {L}egendre function of the
  first kind with respect to its degree''.
\newblock {\em J. Phys. A}, 40(49):14887--14891, 2007.

\bibitem{Szmy2}
R.~{Szmytkowski}.
\newblock {A note on parameter derivatives of classical orthogonal
  polynomials}.
\newblock {\em ArXiv e-prints}, January 2009.

\bibitem{Szmy3}
R.~{Szmytkowski}.
\newblock {On the derivative of the associated Legendre function of the first
  kind of integer degree with respect to its order (with applications to the
  construction of the associated Legendre function of the second kind of
  integer degree and order)}.
\newblock {\em {Journal of Mathematical Chemistry}}, {46}({1}):{231--260},
  {June} {2009}.

\bibitem{Szmy1}
R.~{Szmytkowski}.
\newblock {On the derivative of the associated Legendre function of the first
  kind of integer order with respect to its degree}.
\newblock {\em ArXiv e-prints}, July 2009.

\bibitem{Whip}
F.~J.~W. Whipple.
\newblock {A symmetrical relation between Legendre's functions with parameters
  $\cosh\alpha$ and $\coth\alpha$}.
\newblock {\em {Proceedings of the London Mathematical Society}},
  {16}:{301--314}, {1917}.

\end{thebibliography}

{\it Note added in proof.}---  It has recently come to our attention that the argument 
domain of applicability of the formulas for the derivatives of the Legendre functions 
that we have presented in this paper are actually valid in the complex domain
$z\in{\mathbb C}\setminus[-1,1]\supset \{z: |z|>1, z\in{\mathbb C}\}$.

\end{document}